\documentstyle[amscd,12pt]{amsart}

\newtheorem{thm}{Theorem}[section]

\newtheorem{conj}[thm]{Conjecture}
\theoremstyle{remark}
\newtheorem{rem}[thm]{Remark}

\theoremstyle{definition}

\newcommand{\Hom}{\mathop{\operatorname{Hom}}\nolimits}

\newcommand{\Ker}{\mathop{\operatorname{Ker}}\nolimits}
\newcommand{\im}{\mathop{\operatorname{Im}}\nolimits}
\newcommand{\ad}{\mathop{\operatorname{ad}}\nolimits}
\newcommand{\tr}{\mathop{\operatorname{Tr}}\nolimits}

\newcommand{\Lie}{\mathop{\operatorname{Lie}}\nolimits}
\newcommand{\sgrad}{\mathop{\operatorname{sgrad}}\nolimits}
\newcommand{\vol}{\mathop{\operatorname{vol}}\nolimits}
\newcommand{\supp}{\mathop{\operatorname{supp}}\nolimits}
\newcommand{\curv}{\mathop{\operatorname{curv}}\nolimits}
\newcommand{\rank}{\mathop{\operatorname{rank}}\nolimits}
\newcommand{\dist}{\mathop{\operatorname{dist}}\nolimits}

\begin{document}
\pagestyle{empty}

\title{Jeffrey-Kirwan-Witten Localization Formula  for Reductions at
Regular Co-adjoint Orbits}
\author{Do Ngoc Diep}

\address{Institute of Mathematics, National Centre for Natural Science and 
Technology, P. O. Box 631, Bo Ho, VN-10.000 Hanoi, Vietnam}
\curraddr{International Centre for Theoretical Physics, ICTP P. O. Box
586, 34100, Trieste, Italy}
\email{dndiep@@thevinh.ncst.ac.vn}
\thanks{The author was supported in part by Vietnam National Research
Programme in Fundamental Sciences and in part by the International Centre 
for Theoretical Physics,} 
\subjclass{Primary 22E41, 19E20; Secondary
57T10} 
\dedicatory{}
\keywords{de Rham cohomology, Chern-Weil homomorphism}

\begin{abstract} 
For Marsden-Weinstein reduction at the point $0$ in ${\mathfrak g}^*$, the
well-known Jeffrey-Kirwan-Witten localization formula was proven and then
by M. Vergne modified. We prove in this paper the same kind formula for
the reduction at regular co-adjoint
orbits by using the universal orbital formula of characters. 
\end{abstract}
\maketitle

\section{Introduction and Statement of Results}

Let $(M,\sigma)$ be a symplectic manifold, $G$ a compact Lie group acting
on $M$ by an Hamiltonian action, ${\mathfrak g} = \Lie(G)$ its Lie algebra
and ${\mathfrak g}^*= \Hom_{\mathbb R}({\mathfrak g},{\mathbb R})$. It is
well-known the co-adjoint action of $G$ on ${\mathfrak g}^*$.  Let us consider
the moment map $\mu : M \to {\mathfrak g}^*$ defined by the following formula
 $$\langle \mu(m),X\rangle = f_X(m),$$ where by definition $f_X$ is the
function such that $\imath(X_M)\sigma = df_X$ and $X_M$ is the Hamiltonian
field defined by the action of Lie group $G$ on $M$ $$X_M(m) :=
\frac{d}{dt}|_{t=0}\exp(-tX)m.$$ Let us consider a co-adjoint orbit
${\mathcal O} \in {\mathfrak g}^*/G$.  Recall that $\mu^{-1}({\mathcal O}) 
\longrightarrow {\mathcal O} \subset {\mathfrak g}^*$ is a principal
bundle with the structural group $G$. The quotient $M^{\mathcal O}_{red}
:= G\setminus \mu^{-1}({\mathcal O})$ is known as the {\it
Marsden-Weinstein reduction} of the symplectic manifold $(M,\sigma)$ at
the orbit ${\mathcal O}$ with respect to the moment map $\mu$.  Let us
consider some 1-form $\mu(X)$ on $M$,defined by
 $$ \mu(X)(m) := \langle \mu(m), X\rangle.$$ One defines also an extended
differential form $\sigma_{\mathfrak g} := \mu(X) + \sigma$ and denotes the
horizontal component of this form on $M^{\mathcal O}_{red}$ by
$\sigma^{\mathcal O}_{red}$. 

Assume that some open tubular neighborhood $M^{\mathcal O}_0$ of the
orbit ${\mathcal O}$ is contained in the set of regular values of the
moment maps. This is equivalent also to the assumption that the action of
$G$ on $M^{\mathcal O}_0$ is free. Let us consider the function
$$\mu_{\mathcal O} := \min_{\lambda\in {\mathcal O}} \Vert \mu(X) -
\lambda \Vert_{{\mathfrak g}^*} = \dist_{{\mathfrak g}^*}(\mu(X),
{\mathcal O}).$$
Following Witten, we consider also the form $\frac{1}{2}
\Vert\mu_{\mathcal O}\Vert^2$, which is a $G$-invariant function. Choose
on $M$ an
$G$-invariant metrics $(.,.)$. Denote by $H_{\mathcal O}$ the corresponding Hamiltonian
field, i.e. the symplectic gradient $\sgrad(\frac{1}{2}\Vert\mu_{\mathcal O}\Vert^2)$
and by $\lambda^{M_{\mathcal O}}(.) := (H_{\mathcal O},.)$ the
corresponding $G$-equivariant
differential 1-form. If $\alpha$ is some $G$-equivariant differential
form, let us denote its horizontal component by $\alpha^{\mathcal O}_{red}$ and refer
to it as its reduction on $M^{\mathcal O}_{red}$. Consider $G$-equivariant
differential form of type $e^{i\sigma_{\mathfrak g}(X)}\beta(X),$ where
$\beta(X)$ is closed $G$-invariant differential form depending
polynomially on $X\in{\mathfrak g}$, $\alpha^{\mathcal O}_{red} =
e^{i\sigma^{\mathcal O}_{red}}\beta^{\mathcal O}_{red}.$ Denote $d_X = d -
\imath(X_M)$ the
differential on $G$-equivariant differential forms. We decompose the
manifold $M$ into a union $M = M^{\mathcal O}_0 \cup (M\setminus
M^{\mathcal O}_0)$. 

For each $t\in {\mathbb R}$, $X\in{\mathfrak g}$, consider the generalized
function $\Theta(M,t)$ given by the integral $$\Theta(M,t)(X) :=
\int_Me^{-itd_X\lambda^{M_{\mathcal O}}}\alpha(X).$$ 

\begin{thm}[Main Result] For every
closed $G$-equivariant differential form $\alpha$, there exist limits in
sense of the theory of generalized functions $$\Theta_0^{\mathcal O} :=
\lim_{t\to\infty} \Theta(M_0^{\mathcal O},t),$$
 $$\Theta_{out}^{\mathcal O} := \lim_{t\to\infty} \Theta(M\setminus
M_0^{\mathcal O},t)$$ and one has a decomposition of the integral $\int_M
\alpha$ into their sum $$\int_M \alpha = \Theta_0^{\mathcal O} +
\Theta_{out}^{\mathcal O}$$ also in the sense of generalized functions. 
For each test function $\Phi$, the first summand can be written as
$$\int_{\mathfrak g} \Theta_0^{\mathcal O}(X) \Phi(X) dX = (2\pi i)^{\dim
G}\int_{P^{\mathcal O}} \alpha^{\mathcal O}_{red} \Phi(\Omega)\wedge\vol_\omega.$$
This term can be also expressed by the Kirillov orbital formula for characters
$$\int_{\mathfrak g}\Theta_0^{\mathcal O}(X)\Phi(X)dX = $$
$$(2\pi i)^n \int_{M_{red}^{\mathcal O}}\alpha^{\mathcal O}_{red}\int_{{\mathfrak
g}^*/G} dP({\mathcal O})\int_{\mathcal
O}e^{-i(\Omega,\xi)}(\int_{\mathfrak
g} e^{i\langle \xi,X\rangle}J^{-1/2}_{\mathfrak
g}(X)\Phi(X)dX)d\beta_{\mathcal
O}(\xi)\vol_\omega,$$ 
where $\beta_{\mathcal O}(\xi)$ is the Liouville measure on the co-adjoint orbit,
$J^{-1/2}_{\mathfrak g}(X):= \left(\det{\left( \frac{\sinh(\ad X/2)}{\ad X/2} 
\right)}\right)$ is the Kirillov factor for its universal
character formula and $dP({\mathcal O})$ is the Plancher\`el measure.  
\end{thm}

\begin{rem} The expression of the double integral inside the last formula is
just the Kirillov universal formula for characters of unitary
representations corresponding to the co-adjoint orbit ${\mathcal O}$.  
\end{rem}

\begin{rem} The outer term is some integral over $M_{out}^{\mathcal O}$
and is given by integral also $$\Theta^{\mathcal O}_{out}(X) =
\int_{M^{\mathcal O}\cup C^{\mathcal O}}e^{-i d_X \widetilde{\lambda^{M_{\mathcal
O}}}}
\alpha(X).$$
\end{rem}

These results play some important role in the Witten intersection theory
of cohomology of $G$-equivariant differential forms. 

Witten conjectured the following formula

\begin{conj} $$\int_{\mathfrak g}\Theta_0^{\mathcal O}(X)\Phi(X)dX = (2\pi
i)^{\dim G}\vol(G)\int_{M^{\mathcal O}_{red}}\alpha^{\mathcal O}_{red}W(\Phi),$$
where
$W : C^{\infty}({\mathfrak g})^G \to H^*(M^{\mathcal O}_{red})$ is the
Chern-Weil homomorphism, associated to the principal fibration
$\mu^{-1}({\mathcal O}) \to M^{\mathcal O}_{red}$.  \end{conj}

\begin{rem} The localization formula and, in particular, this conjecture
were proved for $G={\mathbb S}^1$ by Kalkman \cite{kalkman} and Wu
\cite{wu}, for general $G$ and ${\mathcal O}=\{0\}$ by Jeffrey, Kirwan
\cite{jeffreykirwan} and then modified by M. Vergne \cite{vergne}. We
prove this formula by using the universal orbital formula for characters
by Kirillov \cite{kirillov}.  \end{rem}

\section{Proof of Theorem 1.1}

\subsection{Local Fourier Transform}

Let us recall the main moments in Vergne's modification for
Jeffrey-Kirwan-Witten localization theorem. Recall that the group $G$ acts
on $M$ by a Hamiltonian action. For $X\in{\mathfrak g}$, $\xi\in
{\mathfrak g}^*$ and $\Phi\in C^{\infty}({\mathfrak g})$, one has the
Fourier transform $${\mathcal F}(\Phi)(\xi) := (2\pi)^{\dim
G}\int_{\mathfrak g}\Phi(X) e^{-i\langle \xi, X\rangle}dX$$ and by the
well-known inverse Fourier transform, we have $$\Phi(X) = \int_{{\mathfrak
g}^*} e^{i\langle\xi, X\rangle}{\mathcal F}(\Phi)(\xi) d\xi.$$ We can
identify each element $P$ from the symmetric algebra $S({\mathfrak g}^*)$
with a polynomial function on ${\mathfrak g}$, $X \to P(X)$ and also with
a differential operator $P(\partial_{\xi})$ on ${\mathfrak g}^*$, defined
by the property $$P(\partial_{\xi})(e^{\langle\xi,X\rangle}) = P(X)
e^{\langle\xi,X\rangle}.$$ Similarly, we can identify the symmetric
algebra $S({\mathfrak g})$, either with the algebra of polynomial
functions on ${\mathfrak g}^*$ or with the algebra of differential
operators on ${\mathfrak g}$. In particular, to $X\in {\mathfrak g}$
corresponds the Hamiltonian vector field $X_M$ on  $M$.

Consider the algebra ${\mathcal A}^{\infty}_G({\mathfrak g}, M) :=
C^{\infty}({\mathfrak g},{\mathcal A}(M))^G$ of smooth $G$-equivariant
functions $\alpha$, $X\in {\mathfrak g}\mapsto \alpha(X)$. We refer to
${\mathcal A}^{\infty}_G(M)$ as the space of smooth $G$-equivariant
differential forms on $M$, i.e.  $$\alpha(g.X)(g.m) \equiv \alpha(X)(m),
\forall X\in {\mathfrak g}.$$ In particular, our moment maps $\mu : M \to
{\mathfrak g}^*$ defines the function $\mu(X)$, $\mu(X)(m) := f_X(m)$,
$X\in {\mathfrak g}$ as an element of ${\mathcal A}({\mathfrak g},
C^{\infty}(M))$, $C^{\infty}(M) = {\mathcal A}^0(M)$. 

One defines $G$-equivariant co-boundary operator $$d_{\mathfrak g} : 
{\mathcal A}^{\infty}_G({\mathfrak g},M) \to {\mathcal A}^{\infty
+1}_G({\mathfrak g},M)$$ by the formula $$d_{\mathfrak g}\alpha(X) =
d(\alpha(X)) - \imath(X_M)\alpha(X),$$ where the second term is the
contraction of Hamiltonian vector field $X_M$ with the differential form
$\alpha(X)$. It is easy to check that $$d_{\mathfrak g}\alpha(X) =
d^2_{\mathfrak g}\alpha(X) - \imath(X_M)d\alpha(X) - d\imath(X_M)\alpha(X) 
+ \imath(X_M)\imath(X_M)\alpha(X).$$ The first and the last terms are 0 in
virtue of properties of differential forms $\alpha(X)$. The sum of two
remain terms gives us the Lie covariant derivative $ \Lie_{X_M}\alpha(X)$
of $\alpha(X)$ along the vector field $X_M$, which vanishes because
$\alpha$ is $G$-equivariant. 

One denotes also $$d_X := d - \imath(X_M).$$ One has therefore a complex
$G$-equivariant differential forms $$\CD 0 @>>> {\mathbb R} @>>> {\mathcal
A}^0_G(M) @>d_{\mathfrak g} >> {\mathcal A}^1_G(M) @>d_{\mathfrak g}>>
\ldots \endCD$$ It is easy to see that $d_{\mathfrak g}^2 = 0$ and one
defines the $G$-equivariant cohomology as the cohomology of this complex
$${\mathcal H}^{\infty}_G({\mathfrak g},M) = \Ker d_{\mathfrak g} / \im
d_{\mathfrak g}.$$

Suppose that $M$ is oriented and $\alpha\in {\mathcal
A}^{\infty}_G({\mathfrak g},M)$ such that $\forall X\in {\mathfrak g}$,
$\supp \alpha(X)$ is contained in a compact. One defines then $$\int_M
\alpha(X) = \int_M \alpha(X)_{[n]}.$$ if $n = \dim M$ and $\alpha(X)$ is
decomposed into a sum $$\alpha(X) = \alpha(X)_{[0]} + \alpha(X)_{[1]} +
\dots + \alpha(X)_{[n]}$$ of homogeneous differential forms
$\alpha(X)_{[i]}$ of degree $i$.

This means that one has a map $$\int_M : {\mathcal A}_G({\mathfrak g},M) 
\to C^\infty({\mathfrak g})^G.$$ Suppose that our manifold $M$ is closed,
$\partial M = \emptyset$. Then because $\imath(X_M)\alpha(X)$ is of degree
in $-1$ lower and $d\alpha(X)$ is of degree in $ +1$ upper than the degree
of $\alpha(X)$, we have $$ \int_M d_{\mathfrak g}\alpha(X) =
\int_M(d\alpha(X) - \imath(X_M)\alpha(X)) 
                       = \int_M d\alpha(X) = \int_{\partial M} \alpha(X) =
0.$$ Thus we obtain a well-defined map $$\int_M : {\mathcal
H}^\infty_G({\mathfrak g},M) \to C^\infty({\mathfrak g})^G.$$ Let us
consider the case where $M={\mathfrak g}^*$. In this case we have element
$\xi\in {\mathcal A}^\infty_G({\mathfrak g},{\mathfrak g}^*)$, defined as
$X \mapsto (\xi,X)$ and if $U$ is an open subset of ${\mathfrak g}^*$ and
$\beta\in {\mathcal A}^\infty_G({\mathfrak g},U)$, we can consider the
form $\alpha(X) := e^{i(\xi,X)}\beta(X) $ and its differential
$$(d_{\mathfrak g}\alpha)(X) = e^{i(\xi,X)}(i(d\xi,X)+(d_{\mathfrak
g}\beta)(X))),$$ where $d\xi = \sum d\xi^ie^*_i,$ if $X= \sum x_ie^i$ in a
basis $E_1, \dots,E_n$ of ${\mathfrak g}$ and $\xi = \sum \xi^i E^*_i$ in
the corresponding dual basis basis $E^*_1, \dots,E^*_n$ of ${\mathfrak
g}^*$. If $\beta\in {\mathcal A}^{pol}_G({\mathfrak g},U)$ is a polynomial
function on $X$, then $d_{\mathfrak g}\alpha(X) = e^{i(\xi,X)}\gamma(X)$,
where $\gamma$ polynomially depend on $X$.  We can therefore consider the
sub-complex $${\mathcal A}^{\mathcal F}_G({\mathfrak g},U) :=
\{\alpha(X)=e^{i(\xi,X)}\beta(X); \beta\in {\mathcal A}^{pol}_G({\mathfrak
g},V) \}$$ and its cohomology is denoted as ${\mathcal H}_G^{\mathcal
F}({\mathfrak g},U)$. Choose an orientation on ${\mathfrak g }^*$, we have
$$\int_{{\mathfrak g}^*} \alpha(X) = \int_{{\mathfrak g}^*}
\alpha(X)_{[n]},$$ if for example, $\alpha(X)$ is a rapidly decreasing
$C^\infty$-functions on ${\mathfrak g}^*$. The result is also a rapidly
decreasing $C^\infty$-function on ${\mathfrak g}$. We can therefore define
also its Fourier transform. Suppose that $$\alpha(X)_{[n]} = \sum
P_a(X)\alpha_a(\xi)d\xi,$$ where $Pa \in {\mathbf S}({\mathfrak g}^*)$,
$\alpha_a(\xi) \in C^\infty({\mathfrak g}^*).$ Then we have a formula
$${\mathcal F}(\int_{{\mathfrak g}^*} \alpha) = (\sum
P_a(i\partial_\xi)\alpha_a(\xi))d\xi.$$ M. Vergne defined the local
Fourier transform as the generalized function $$V(\alpha)= (\sum
P_a(i\partial_\xi)\alpha_a(\xi))d\xi.$$ It was proven that if $\beta\in
{\mathcal A}^{\mathcal F}_G({\mathfrak g},V)$, then $V(d_{\mathfrak
g}\beta) = 0$.  The local Fourier transform is defined indeed on the
cohomology groups $$V : {\mathcal H}^{\mathcal F}_G({\mathfrak g},U) \to
{\mathcal A}^n(U)^G.$$

Let us consider the moment map $\mu : M \to {\mathfrak g}^*$. For each
element $m\in M$, $\mu(m)$ is a function on ${\mathfrak g}^*$, $(\mu(m),X) 
= f_X(m)$ and for each $X$, one has a function $(\mu(.),X)$ on $M$. Thus
one can define the function $X \mapsto e^{i(\mu,X)} \in {\mathcal
A}^\infty_G({\mathfrak g},M)$ and for all $\beta\in {\mathcal
A}^{pol}_G({\mathfrak g},M)$, $\alpha(X) = e^{i(\mu,X)}\beta(X) \in
{\mathcal A}^\infty_G({\mathfrak g},M)$. The subspaces $${\mathcal
A}^\mu_G({\mathfrak g},M) :=
 \{\alpha(X) = e^{i(\mu,X)}\beta(X) ; \beta\in {\mathcal
A}^{pol}_G({\mathfrak g},M)  \}$$ is stable under $d_{\mathfrak g}$. We
have therefore a sub-complex $({\mathcal A}^\mu_G({\mathfrak
g},M),d_{\mathfrak g})$ and the corresponding cohomology is denoted
 by ${\mathcal H}^{pol}_G({\mathfrak g},M)$.  \begin{thm}[\cite{vergne}]
Assume that the manifold $M$ is oriented and $U$ is a $G$-invariant open
set containing in the set of regular values of the moment map $\mu$. Let
$\alpha \in {\mathcal A}^{\mu}_G({\mathfrak g},M)$. Then over $U$ one has
$${\mathcal F}(\int_M\alpha) = V(\mu_*\alpha),$$ where $\mu_*\alpha$ is
the push-forward of $\alpha$. If $\alpha$ is closed, ${\mathcal
F}(\int_M\alpha)$ depend only on the cohomology class of $\alpha$ in
${\mathcal H}^\mu_G({\mathfrak g},\mu^{-1}(U))$. Thus for $\alpha\in
{\mathcal H}^{\mu}_G({\mathfrak g},\mu^{-1}(U))$, in order to determine
${\mathcal F}(\int_M \alpha)$ near a regular value $F$ of $\mu$, we need
only to determine the class of $\alpha$ in ${\mathcal
H}^{\mu}_G({\mathfrak g},\mu^{-1}(U))$, where $U$ is a $G$-invariant tub
ular open neighborhood of the orbit ${\mathcal O}$ of $F$.  \end{thm}

Let us consider also push-forward $\mu_*((\omega_a)_{[\dim M]})$, which is
a Radon measure. We have by definition, $$\int_M \omega = \int_{{\mathfrak
g}^*} \mu_*\omega .$$ One has also $$(\int_M\alpha)(X) = \sum_a P_a(X)
\int_{{\mathfrak g}^*}e^{i(\mu,X)}\mu_*((\omega_a)_{[\dim M]})$$ and
$${\mathcal F}(\int_M \alpha) = \sum_a
P_a(i\partial_\xi)(\mu_*(\omega_a)_{[\dim M]}) = V(\mu_*\alpha).$$ From
these, one has a result 

\begin{thm}[Berline-M. Vergne \cite{vergne}, 1982]
Let $T \subset G$ be a maximal torus of $G$ and $M^T$ is the sub-manifold
of fixed point under the action of the torus $T$, then $\int_M\alpha$
depend only on the restriction $\alpha|_{M^T}$.  \end{thm}

Let us now deduce an explicit formula for $V\mu_*$ near the point $F\in
{\mathcal O}\subset {\mathfrak g}^*$ Suppose that the action of $G$ on
$\mu^{-1}({\mathcal O})$ is locally free.  Choose a $G$-invariant Euclidean
norm $\Vert . \Vert$ on ${\mathfrak g}$ and $U$ a $G$-invariant open ball
centered at $F\in {\mathfrak g}^*$, $N^{\mathcal O} := P^{\mathcal O}
\times_G U$ is a fibration over $N^{\mathcal O}_{red} := G\setminus
P^{\mathcal O}$. The map $\mu : N^{\mathcal O} \to U$ is the projection on
the second component. Let us consider a form $\alpha\in {\mathcal
A}^\mu_G({\mathfrak g},N^{\mathcal O}),$ the restriction
$\alpha|_{P^{\mathcal O}}$ is $G$-equivariant. Choose an orientation in
$P^{\mathcal O} = \mu^{-1}({\mathcal O})$ and a basis $E_1, E_2,\dots,E_n$
of ${\mathfrak g}$ and the corresponding dual basis $E^*_1,
E^*_2,\dots,E^*_n$ of ${\mathfrak g}^*$. Choose a connection, i.e. a
trivialization by open covering of the base $N ^{\mathcal O}_G$.  Let us
write $\omega$ a connection on the fibration $P^{\mathcal O} \to
G\setminus P^{\mathcal O}$, $$\omega = \sum \omega_k E_k,$$ and $$\Omega =
\curv(\omega) = \sum \Omega_k E_k$$ the curvature of $\omega$ with values
in $\mathfrak g$. It is not hard to see that $q^*: {\mathcal
H}^\infty_G({\mathfrak g},P^{\mathcal O}) \cong H^*(G\setminus P^{\mathcal
O}). $ If $\Phi$ is a polynomial function on ${\mathfrak g}$ then
$\Phi(\Omega$ is a differential form on $P^{\mathcal O}$. Because
$\alpha(X)$ is $G$-equivariant, then $\alpha(\Omega)$ is a form on
$P^{\mathcal O}$ and if $\alpha$ is closed $G$-equivariant, the horizontal
component $\alpha_{red} := h(\alpha(\Omega))$ is a closed differential
form on $N^{\mathcal O}_{red}$. The map $\alpha \mapsto \alpha_{red}$ is
just the inverse map $$(q^*)^{-1} : {\mathcal H}^\infty_G({\mathfrak
g},P^{\mathcal O}) \to H^*_{DR}(G\setminus P^{\mathcal O}).$$ one denote
also $(\Omega,\xi) = \sum \Omega^k \xi_k$ on $P^{\mathcal O}$ and
$\vol^{\mathcal O}_\omega = \omega_1 \wedge \dots\wedge \omega_n$ be the
vertical form on $P^{\mathcal O}$ of degree $n =\dim G$. Let us denote
also by coordinates $$dX = dx_1\wedge dx_2 \wedge \dots \wedge dx_n$$ the
volume form then $E_1^* \wedge E_1 \wedge \dots \wedge E_n^* \wedge E_n =
dX\wedge d\xi$, where $d\xi = (-1)^{n(n+1)/2}d\xi^1\wedge d\xi^2\dots
\wedge d\xi^n$. 

\begin{thm}[\cite{vergne}] If $P^{\mathcal O}$ is oriented with free
$G$-action then for all $\Phi\in{\mathbf S}({\mathfrak g}^*)^G,$ one has
$$((V\mu_*\alpha)/d\xi, \Phi) = i^n\int_{P^{\mathcal
O}}\alpha^{\mathcal O}_{red}\Phi(-i\Omega)\vol_\omega= i^n\vol(G)\int_{N^{\mathcal
O}_{red}}\alpha^{\mathcal O}_{red}\Phi(-i\Omega).$$ \end{thm}

\begin{thm}[Jeffrey-Kirwan\cite{jeffreykirwan}]\label{jeffreykirwanformula} For
$\alpha\in {\mathcal
A}^\mu_G({\mathfrak g},M),$ $${\mathcal F}(\int_M\alpha) =
i^n\int_{P^{\mathcal O}} \alpha_{red} e^{-i(\xi,\Omega)}\vol_\omega.$$
\end{thm} Witten defined the integral $$Z(\varepsilon) =
\int_M\int_{\mathfrak g} e^{i\sigma_{\mathfrak g}(X)} \beta(X)
e^{-\frac{\varepsilon\Vert X\Vert^2}{2}}dX.$$

\begin{thm}[Witten\cite{vergne},\cite{jeffreykirwan}] Let
$(M,\sigma,\mu)$be a symplectic manifold with a Hamiltonian action of $G$,
${\mathcal O}$ a co-adjoint orbit of $G$ in ${\mathfrak g}^*$,
$\mu_{\mathcal O} : \mu^{-1}({\mathcal O}) \to M^{\mathcal O}_{red} = G
\setminus \mu^{-1}({\mathcal O})$, $R$ the smallest critical value of
$\Vert\mu_{\mathcal O}\Vert^2$, $r < R$. Then for each closed
$G$-equivariant polynomially depend on $X$, there exists a constant $C$
such that $$Z(\varepsilon) = (2\pi i)^{\dim G} \int_{M^{\mathcal O}_{red}}
e^{i\sigma_{red}} \beta_{red} e^{\frac{\varepsilon \vert \Omega\vert^2}{2}}
+ N(\varepsilon),$$ where $$\vert N(\varepsilon \vert \leq C
e^{-r/2\varepsilon},\forall \varepsilon > 0.$$ 
\end{thm}

\subsection{Kirillov Universal Character Formula}

Let us describe in this section the well-known Kirillov orbital formula
for characters of representations.  

\begin{thm}[Liouville Measure]
$$\int_{M_F}\Phi d\beta_F = \prod_{\alpha\in P_F} \frac{\langle
F,iH_\alpha\rangle}{2\pi }\int _{G/G_F}\Phi(gF) d\bar{g},$$ where
$d\bar{g}$ is the quasi-invariant measure on the co-adjoint orbit
${\mathcal
O} = G_F \setminus G$. 
\end{thm}

The theorem what follows describes the image of the Fourier transform of
Liouville measure.

\begin{thm}[Harish-Chandra, see \cite{berlinegetzlervergne}, Theorem 7.24
and Corollary 7.25]
Given $\lambda \in {\mathfrak t}^*$, let $W_\lambda = \{ w\in W \vert
w\lambda = \lambda \}$ be the stabilizer of $\lambda$ in the Weyl group
$W$. For $X\in {\mathfrak t}$, $X$ regular, the Fourier transform
$$F_{M_\lambda}(X) = \int_{M_\lambda} e^{if(X)}d\beta_\lambda(f)$$ is given by
the formula $$F_{M_\lambda}(X) = \sum_{W/W_\lambda} \frac{e^{i\langle
w\lambda, X\rangle}}{\prod_{\alpha\in P_\lambda}\langle
w\alpha,X\rangle}.$$

If $M_\lambda$ is a regular orbit of the co-adjoint representation, then
for $X\in {\mathfrak t},$ $X$ regular,
$$F_{M_\lambda}(X) = \prod_{\alpha\in P_\lambda} \langle \alpha, X\rangle
^{-1} \sum_{w\in W} \varepsilon(w) e^{i\langle w\lambda, X\rangle}.$$
\end{thm}

Let us recall that if $G$ is a connected compact Lie group, we denote by $T \subset
G$ a fixed maximal torus in $G$, $$L_T := \{ X\in {\mathfrak t} \vert e^X = 1 \}$$
the lattice of co-roots, $$L_T^* = \{\lambda \in {\mathfrak t}^* \vert \lambda(X) \in
2\pi {\mathbb Z}, \forall X \in L_T \} \subset {\mathfrak t}^*$$ the root lattice.
One fixes some basis in this lattice and every element from lattice can be expressed
as linear combination with integral coefficients. One fixes a system $P$ of positive
roots in $\Delta$. Let us denote by $\rho$ the half-sum of positive roots and
$$X_G
:= i\rho + L_T^*.$$ E. Cartan described the irreducible finite dimensional
representations of $G$ and H. Weyl described the characters of irreducible
representations: There is a bijective correspondence between regular elements
$\lambda \in X_G$ and the irreducible representations, denoted by $T_\lambda$.
Weyl character formula is
$$\tr(T_\lambda(e^X)) = \frac{\sum_{w\in W} \varepsilon(w) e^{i\langle
w\lambda,X\rangle}}{\prod_{\alpha\in P_\lambda}(e^{\langle
\alpha,X\rangle/2}-e^{-\langle \alpha,X\rangle}/2)},$$

Let us consider the function
$$J_{\mathfrak g}(X) := \det \left(\frac{\sinh(\ad X)/2}{\ad X/2}
\right).$$ Then
for all $X \in {\mathfrak t},$
$$J_{\mathfrak g}(X) = \prod_{\alpha\in\Delta}
\frac{e^{\langle\alpha,X\rangle/2}-e^{-\langle\alpha,X\rangle/2}}{\langle\alpha,X\rangle}$$
$$ = (\prod_{\alpha\in\Delta}
\frac{e^{\langle\alpha,X\rangle/2}-e^{-\langle\alpha,X\rangle/2}}{\langle\alpha,X\rangle})^2.$$

\begin{thm}[Chevalley's Theorem, see \cite{berlinegetzlervergne}, Thm. 7.28]
There are isomorphisms between $G$-invariants and $W$-invariants:
$$C^\infty({\mathfrak t})^W \cong C^\infty({\mathfrak g})^G,$$
$$C^\omega({\mathfrak t})^W \cong C^\omega({\mathfrak g})^G,$$
$$C[{\mathfrak t}]^W \cong C[{\mathfrak g}]^G.$$
\end{thm}

The main reason for these isomorphisms is the correspondence $\Phi \mapsto c(\Phi)$,
where 
$$c(\Phi) = \frac{1}{\vert W\vert}\int_G (\partial_w(\pi_{{\mathfrak
g}/{\mathfrak t}}\Phi)(pr(g.X))dg,$$ where
$\partial_w := \prod_{\alpha\in P}\langle \rho,\alpha\rangle^{-1}\prod_{\alpha\in P} 
\partial_{iH_\alpha}$, $\pi_{{\mathfrak g}/{\mathfrak t}}(X) := \prod_{\alpha\in P} 
i\alpha(X) \in {\mathbb C}[{\mathfrak t}]$.   

From this theorem, it is easy to see that there is a unique function 
$$J^{1/2}_{\mathfrak g}(X) := \prod_{\alpha\in P} \frac{\sinh\langle
\alpha,X\rangle/2}{\langle\alpha,X\rangle/2}$$ such that
$$J^{1/2}_{\mathfrak g}(X)\tr(T_\lambda(e^X)) = \sum \langle\alpha,X\rangle^{-1}
\sum_{w\in W} \varepsilon(w) e^{i\langle w\lambda,X\rangle}$$
$$ =  F_{M_\lambda}(X) = \int_{M_\lambda} e^{if(X)}d\beta_\lambda(f).$$

The same formula for semi-simple Lie groups was obtained by Rossman:  

\begin{thm}[see \cite{berlinegetzlervergne}, Theorem 7.29] Let $M$ be a closed
co-adjoint orbit of a real semi-simple Lie group $G$ with non-empty discrete
series, i.e. $\rank(G) = \rank(K)$ and let $W = W({\mathfrak k}_{\mathbb C},
{\mathfrak t}_{\mathbb C})$ be the compact Weyl group. Then for regular element $X\in
{\mathfrak t}_r$, we have the following results:
\begin{enumerate}
\item If $M \cap {\mathfrak t}^* = \emptyset$, then $F_M(X) = 0$.
\item If $M = G.\lambda$ with $\lambda \in {\mathfrak t}^*$, and $W_\lambda$ is the
subgroup of $W$ stabilizing $\lambda$, then
$$F_{M_\lambda}(X) = (-1)^{n(\lambda)} \sum_{W/W_\lambda} \frac{e^{i\langle
w\lambda,X\rangle}}{\prod_{\alpha\in P_\lambda}\langle w\alpha,X\rangle}.$$
\end{enumerate}
\end{thm}

\subsection{Witten Type Localization Theorem}

Let $(M,\sigma,\mu)$ be a compact symplectic manifold with a Hamiltonian action of
compact Lie group $G$. Let us assume that the co-adjoint orbit ${\mathcal O}$ is
contained in the set of regular values of $\mu$. Assume that the action of $G$ is
free in $P^{\mathcal O} = \mu^{_1}({\mathcal O})$. Let us denote $\omega$ a
connection form on $P^{\mathcal O}$ with curvature $\Omega$. We consider the so
called Marsden-Weinstein reduction $M^{\mathcal O}_{red} := G \setminus P^{\mathcal
O}$ of $M$. We denote $\alpha^{\mathcal O}_{red}$ the de Rham cohomology class
$(\alpha\vert_{P^{\mathcal O}})_{red}$ on $M^{\mathcal O}_{red}$ determined by
$\alpha\vert_{P^{\mathcal O}}$. In particular $(\sigma^{\mathcal O}_{\mathfrak
g})_{red}$ is the symplectic form $\sigma^{\mathcal O}_{red}$ on $M^{\mathcal
O}_{red}$. 

Let us consider the function $\mu_{\mathcal O} := \min_{\lambda \in {\mathcal O}}
\Vert \mu(x) - \lambda \Vert$, the distance from $\mu(X)$ to the co-adjoint orbit
${\mathcal O}$. Following Witten, we introduce the function $\frac{1}{2}\Vert
\mu_{\mathcal O}\Vert^2$ and its Hamiltonian field $H_{\mathcal O}$. This is an
invariant vector field on $M_{\mathcal O}$. Choose a $G$-invariant metric $(.,.)$ on
$M_{\mathcal O}$ and put
$\lambda^{M_{\mathcal O}} := (H_{\mathcal O},.)$. Then $\lambda^M_{\mathcal O}$ is a
$G$-invariant
1-from on $M_{\mathcal O}$. Let $R$ be the smallest critical value of the function $\Vert
\mu_{\mathcal O}\Vert^2$. Let $r < R$ and let
$$M_0^{\mathcal O} = \{ x\in M; \Vert \mu_{\mathcal O}(x)\Vert^2 < r \}, \quad
M_{out}^{\mathcal O} = \{ x\in M; \Vert \mu_{\mathcal O}(x)\Vert^2 > r \}.$$ 

Let $\alpha(X)$ be a closed $G$-invariant differential form on $M$. Let us consider
$$\Theta(M,t) = \int_M e^{_itd_X\lambda^M} \alpha(X),$$ which is independent of $t$
in virtue that $\alpha$ is a closed form and $e^{_itd_X\lambda^M}$ congruent to 1 in
cohomology. Let us consider the integrals
 $$\Theta(M_0^{\mathcal O},t) = \int_{M^{\mathcal O}_0} e^{_itd_X\lambda^M}
\alpha(X),$$
$$\Theta(M_{out}^{\mathcal O},t) = \int_{M_{out}^{\mathcal O}} e^{_itd_X\lambda^M}
\alpha(X).$$
Then the values $\Theta(M_0^{\mathcal O},t)(X)$ and $\Theta(M_{out}^{\mathcal
O},t)(X)$ $C^{\infty}$-smoothly depend on $X\in {\mathfrak g}$. 

\begin{thm}
For every $t\in {\mathbb R}$ and $X\in {\mathfrak g}$, there is a decomposition 
$$\left(\int_M e^{-itd_X\lambda^M}\alpha \right)(X) = \Theta(M^{\mathcal O}_0,t)(X)
+
\Theta(M^{\mathcal O}_{out},t)(X).$$ There exist the limits $\Theta_0^{\mathcal O} =
\lim_{t\to\infty}\Theta(M^{\mathcal O}_0,t)$ and $\Theta_{out}^{\mathcal O} =
\lim_{t\to\infty}\Theta(M^{\mathcal
O}_{out},t)$ in the sense of distributions such that
$$\int_{\mathfrak g}\left(\int_M \alpha \right)(X)\Phi(X) dX =
\int_{\mathfrak g}\Theta_0^{\mathcal O}(X)\Phi(X)
dX +
\int_{\mathfrak g}\Theta^{\mathcal O}_{out}(X)\Phi(X) dX,$$ for every test function
$\Phi$, where
$$\int_{\mathfrak g} \Theta_0^{\mathcal O}(X) \Phi(X) dX = (2\pi i)^{\dim
G}\int_{P^{\mathcal O}} \alpha^{\mathcal O}_{red} \Phi(\Omega)\wedge\vol_\omega.$$ 
\end{thm}

\begin{pf}
The proof is the same as in the note of M. Vergne \cite{vergne} with change
everywhere $M_0$, $M_{out}, \dots, \mu, \dots $  by $M_0^{\mathcal O}$,
$M_{out}^{\mathcal O}, 
\dots, \mu^{\mathcal O}, \dots $
We need : 
\begin{itemize}
\item to prove the existence of limit $\Theta_0 = \lim_{t\to 0} \Theta(M_0,t)$ in
the sense of generalized functions and
\item  to compute it.
\end{itemize}
Let us start now to do this. Let $E_1, \dots, E_n$ be a basis of ${\mathfrak g}$,
such that
$E_1, \dots, E_k$ is a basis of a stabilizer ${\mathfrak g}_F = {\mathcal O}^\perp$
in ${\mathfrak g}$. This means that 
$$\mu_{\mathcal O} = \sum_{i=k}^n \mu(E_i)E_i$$ and 
$$\frac{1}{2}\Vert \mu_{\mathcal O}\Vert^2 = \sum_{i=k}^n \mu(E_i)\mu(E_i) =
\sum_{i=k}^n \mu(E_i)\imath((E_i)_M)\sigma.$$ One deduces easily that 
$$\lambda^{M_{\mathcal O}} = (H_{\mathcal O},.)
 = \sum_{i=1}^m \mu(E_i)\omega^i_M,$$                            
where $\omega^i_M = ((E_i)_M,.)$ and $\omega_M = \sum_{i=1}^m \omega^i_M E_i.$  
We have therefore $$\lambda^{M_{\mathcal O}} = (\omega^{\mathcal O}_M,\mu_{\mathcal
O}).$$ On $M_0^{\mathcal O}$ the action of $G$ is locally free. Thus we an choose on
$M_0^{\mathcal O}$ a metric $(.,.)$ such that $((E_i)_M,(E_j)_M) = \delta_{ij}$ on
$M_0^{\mathcal O}$. One has hence $$\omega_M(X_M) = X, \forall X \in {\mathfrak
g}.$$ This means that $\omega_M$ is a connection form on $M_0^{\mathcal O}$ and 
$$\lambda^{M_{\mathcal O}}(X_M) = (H_{\mathcal O},X_M) = \mu_{\mathcal O}(X).$$
Let us consider the function $f_{\lambda^{M_{\mathcal O}}} : M \to {\mathfrak g}^*$,
defined by the condition that
$$f_{\lambda^{M_{\mathcal O}}}(X) = \lambda^{M_{\mathcal O}}(X_M).$$ Then it is easy
to see that
$$\langle f_{\lambda^{M_{\mathcal O}}},\mu_{\mathcal O}\rangle = \sum_{i=1}^m
\mu_{\mathcal O}(E_i) ((E_i)_M ,H_{\mathcal O}) = (H_{\mathcal O}, H_{\mathcal O})
\geq 0.$$ On $M_0^{\mathcal O}$ we have
$$d_X\lambda^{M_{\mathcal O}} = -i(\mu_{\mathcal O},X) + d\lambda^{M_{\mathcal O}} =
-i(\mu_{\mathcal O},X) + (\mu,d\omega_M) - (\omega_M, d\mu).$$ Thus we have
$$\int_{\mathfrak g} \Theta_0^{\mathcal O}(X)\Phi(X) dX = \int_{M_0^{\mathcal
O}}\int_{\mathfrak g} e^{it(\mu,X) -it(\mu,d\omega_M) +
it(\omega_M,d\mu)}\alpha(X)\Phi(X) dX.$$
Let $\varepsilon > 0$ be a small number. Let us consider a small tubular
neighborhood $M_\varepsilon^{\mathcal O}$ of the orbit ${\mathcal O}$,
$$M^{\mathcal O}_0 = \{ x\in M^{\mathcal O} ; \Vert \mu_{\mathcal O}(x) \Vert <
\varepsilon \}$$ and let $$\chi_{\mathcal O}(m) := \left\{ \begin{array}{ll}
1 & \quad \forall m\in M^{\mathcal O}_{\varepsilon/2}\\  0 &\quad \forall m\notin
M^{\mathcal O}_{\varepsilon/2}\end{array}\right.$$ to be the cut-off function. Then 
$$\lim_{\varepsilon \to 0} \int_{M_0^{\mathcal O}}(1-\chi_{\mathcal
O}(m))(\int_{\mathfrak g} e^{itd_X\lambda^{M_{\mathcal O}}}\alpha(X)\Phi(X)) =0.$$
Let us denote $\beta(X) := \alpha(X) \Phi(X)\in C^\infty_{cpt}({\mathfrak g},
{\mathcal A}(M_{\mathcal O}))$. Then we have
$$\int_{\mathfrak g} e^{-itd_X\lambda^{M_{\mathcal O}}}\alpha(X)\Phi(X) =
\int_{\mathfrak g} e^{it(\mu,X)} e^{-itd\lambda^{M_{\mathcal O}}}\alpha(X)\Phi(X) dX
= e^{-itd\lambda^{M_{\mathcal O}}}\widehat{\beta}(t\mu_{\mathcal O}),$$ where
$\widehat{\beta}$ is, by definition, the Fourier transform of $\beta = \alpha.\Phi$. 
Because $\beta$ is of compact support, its Fourier transform $\widehat{\beta}$
is of Schwartz class and because
$e^{-itd\lambda^{M_{\mathcal O}}}$ is polynomial on $t$, we have
$$\lim_{t\to\infty}\int_{M_0^{\mathcal O}}(1-\chi_{\mathcal O}(m)(\int_{\mathfrak
g} e^{-itd\lambda^{M_{\mathcal O}}}\alpha(X)\Phi(X)dX) =0.$$
Because on the support of the function $1-\chi_{\mathcal O}$, one has the estimate
$\Vert \mu_{\mathcal O}(m) \Vert \geq \frac{1}{2}\varepsilon$, then
$$\lim_{t\to\infty} \int_{M_0^{\mathcal O}}((\int_{\mathfrak
g} e^{-itd\lambda^{M_{\mathcal O}}}\alpha(X)\Phi(X)dX) =$$
$$ =\lim_{t \to
\infty}\int_{M^{\mathcal O}_\varepsilon}\chi_{\mathcal O}(m)(\int_{\mathfrak
g} e^{-itd\lambda^{M_{\mathcal O}}}\alpha(X)\Phi(X)dX).$$
Choose $N^{\mathcal O} = P^{\mathcal O} \times {\mathfrak g}^*$ and let
$\omega_{\mathcal O} = \omega_M^{\mathcal O}|_{P^{\mathcal O}}$, then
$\omega_{\mathcal O}$ is a connection form on $P^{\mathcal O}$. Choose $\varepsilon$
sufficiently small, then $M^{\mathcal O}_\varepsilon \cong \mbox{ open set of }
N^{\mathcal O} = P^{\mathcal O} \times {\mathfrak g}^*$. Thus the map
$\mu_{\mathcal
O} : N^{\mathcal O} \to {\mathfrak g}^*$ becomes the projection $(x,\xi) \mapsto
\xi$. This isomorphism is identity on $P^{\mathcal O}$. Because $\chi_{\mathcal O}$
has compact support contained in $M^{\mathcal O}_\varepsilon$, we can consider the
integral $$\int_{M^{\mathcal O}_\varepsilon}\chi_{\mathcal O}(m)(\int_{\mathfrak
g} e^{-itd\lambda^{M_{\mathcal O}}}\alpha(X)\Phi(X)dX)$$ as an integral over
$N^{\mathcal O}$ and still denote $\omega_{\mathcal O}^M$ for the 1-form
corresponding to $\omega_{\mathcal O}^M$ and $\omega_{\mathcal O}^M|_{P^{\mathcal 
O}} = \omega_{\mathcal O}$.  We have thus
$$\lim_{t\to \infty}\int_{M_0^{\mathcal O}}(\int_{M^{\mathcal O}_\varepsilon} 
\chi_{\mathcal O}(m)(\int_{\mathfrak
g} e^{-itd\lambda^{M_{\mathcal O}}}\alpha(X)\Phi(X)dX) = $$
$$= \lim_{t\to\infty}\chi_{\mathcal O}(m)(\int_{\mathfrak g}
e^{it(\xi,X)}e^{it(\omega_{\mathcal O}^M,d\xi)-it(\xi,d\omega_{\mathcal
O}^M)}\alpha(X) \Phi(X)dX).$$ 
By the same reason as above, there exist polynomials $P_k(t\xi,td\xi)$ in
$\xi^1,\dots,\xi^m$, $d\xi^1, \dots, d\xi^m$ such that 
$$e^{it(\omega^M_{\mathcal O},d\xi) -it(\xi,d\omega^M_{\mathcal O}) = \sum_{k=1}^m 
P_k(t\xi,td\xi)\mu_k},$$ where $\mu_k$ are differential forms on $N^{\mathcal O}$,
independent of $t$. Denote $$\nu_k := \chi_{\mathcal O}\mu_k \wedge
\alpha(X)\Phi(X),$$ we have to study the integrals of type
$$\int_{N^{\mathcal O}}(\int_{\mathfrak g} e^{it(\xi,X)} P_k(t\xi,td\xi)\nu_k(X)
dX).$$ We denote $\nu_0 = \nu(X)|_{P^{\mathcal O}}$, then the map $X \mapsto
\nu_0(X)$ is an element in the class $C^\infty_0({\mathfrak g}, {\mathcal
A}(P^{\mathcal O}))$.
Its Fourier transform $\widehat{\nu}_0(\xi) = \widehat{\nu_0(X)}(\xi)$ as
differential forms
on $N^{\mathcal O} = P^{\mathcal O} \times {\mathfrak g}^*$. It was shown that if
$G(\xi,d\xi)$ is a polynomial, then for all $\nu\in C^\infty_0({\mathfrak g},
{\mathcal A}_{cpt}(N^{\mathcal O}))$, 
$$\lim_{t\to\infty} \int_{N^{\mathcal O}}(\int_{\mathfrak g}
e^{it(\xi,X)}G(t\xi,td\xi))\nu(X)dX = \int_{N^{\mathcal O}}
G(\xi,d\xi)\widehat{\nu}_0(\xi).$$
Because $\chi_{\mathcal O} |_{P^{\mathcal O}} \equiv 0$, $\omega^M_{\mathcal
O}|_{P^{\mathcal O}} = \omega_{\mathcal O}$, we have
$$\lim_{t\to\infty} \int_{N^{\mathcal O}}\chi_{\mathcal O}(m)(\int_{\mathfrak
g}e^{it(\xi,X)}e^{it(\omega^M_{\mathcal O},d\xi)-it(\xi,d\omega^M_{\mathcal
O})}\alpha(X)\Phi(X)dX)=$$
$$= \int_{N^{\mathcal O}}
e^{i(\omega,d\xi)-i(\xi,d\omega)}\widehat{\Phi\alpha_0}(\xi)= 
\int_{N^{\mathcal O}} \int_{\mathfrak g} e^{-id_X\lambda^{M_{\mathcal
O}}}\alpha_0(X)\Phi(X) dX,$$ see (\cite{berlinegetzlervergne},Lemma 21). We have
seen that $$\int_{\mathfrak g} \Theta_0^{\mathcal O}(X)\Phi(X) dX =
\int_{N^{\mathcal O}} \int_{\mathfrak g} e^{-itd_X\lambda^{M_{\mathcal O}}}
\alpha_0(X) \Phi(X) dX.$$
Let us now prove that 
$$\int_{N^{\mathcal O}} \int_{\mathfrak g} e^{-itd_X\lambda^{M_{\mathcal O}}}
\alpha_0(X) \Phi(X) dX = (2\pi i)^{\dim G} \int_{P^{\mathcal O}} \alpha^{\mathcal
O}_{red} \int\Phi(\Omega)\wedge \vol_\omega.$$
Indeed, first we have
$$\int_{N^{\mathcal O}} \int_{\mathfrak g} e^{-itd_X\lambda^{M_{\mathcal O}}}
\alpha_0(X) \Phi(X) dX = $$ $$=\int_{P^{\mathcal O}} \alpha^{\mathcal
O}_{red}(\int_{N^{\mathcal O}/P^{\mathcal O}}\int_{\mathfrak g}
e^{-id_X\lambda^{M_{\mathcal O}}}\alpha_0(X)\Phi(X)dX),$$
with $$e^{-id_X\lambda^{M_{\mathcal 
O}}} = e^{i(\xi,X)}e^{-i(d\omega,\xi)+ i(\omega,d\xi)}.$$ Remark that its
term of maximal  degree in $d\xi$ is
$cd\xi_1\dots  d\xi_m \wedge \vol_\omega$, where $c=i^n\varepsilon$ and
$\varepsilon = (-1)^{n(n+1)/2}$ is a sign. Recall also
that $\Omega = d\omega + \frac{1}{2}[\omega,\omega]$ is the curvature of the
connection form $\omega = \omega_{\mathcal O}^M$. Because $\omega_i \wedge
\vol_\omega = 0, \forall i$, we have
$$e^{-i(d\omega,\xi)}\wedge\vol_\omega = e^{-i(\Omega,\xi)}\wedge \vol_\omega$$ and
then 
$$\int_{N^{\mathcal O}/P^{\mathcal O}} \int_{\mathfrak g}
e^{-itd_X\lambda^{M_{\mathcal O}}}
\alpha_0(X) \Phi(X) dX = $$ $$=c \int_{M^{\mathcal O}/P^{\mathcal O}}
e^{-i(d\omega,\xi)}(\int_{\mathfrak g} e^{i(\xi,X)}\Phi(X) dX)d\xi\vol_\omega $$
$$= c\int_{N^{\mathcal O}/P^{\mathcal O}}
e^{-i(\Omega,\xi)}(\int_{\mathfrak g} e^{i(\xi,X)}\Phi(X)dX)d\xi\vol_\omega.$$
Following the Fourier inverse formula, we have
$$\int_{N^{\mathcal O}/P^{\mathcal O}}
e^{-i(\Omega,\xi)}(\int_{\mathfrak g} e^{i(\xi,X)}\Phi(X)dX) = (2\pi)^n
\Phi(\Omega).$$ The last integral is therefore equal $$(2\pi i)^n \Phi(\Omega)\wedge
\vol_\omega.$$ The proof is therefore achieved.
\end{pf}
\begin{rem}
By the same way as preceded in \cite{vergne}, we can deduce the Jeffrey-Kirwan
formula \ref{jeffreykirwanformula} and the outer term formula. 
\end{rem}

We finish the proof of our main theorem in this section by using together
the local Fourier transform and the universal orbital formula for
characters by Kirillov. 
By the Jeffrey-Kirwan-Witten localization theorem, we have
$$\int_{\mathfrak g} \Theta_0^{\mathcal O}(X) \Phi(X) dX = (2\pi i)^{\dim
G}\int_{P^{\mathcal O}} \alpha^{\mathcal O}_{red} \Phi(\Omega)\wedge\vol_\omega$$
$$=  i^n\varepsilon \int_{N^{\mathcal O}/P^{\mathcal O}}
e^{-i(\Omega,\xi)}(\int_{\mathfrak g} e^{i(\xi,X)}\Phi(X)dX)d\xi\vol_\omega.$$
By the Kirillov universal trace formula, we have
$$ (2\pi)^n \Phi(\Omega) = \int_{N^{\mathcal O}/P^{\mathcal O}}
e^{-i(\Omega,\xi)}(\int_{\mathfrak
g} e^{i(\xi,X)}\Phi(X)dX)  = \int_{{\mathfrak g}^*}
e^{-i(\Omega,\xi)}(\int_{\mathfrak
g} e^{i(\xi,X)}\Phi(X)dX)  $$
$$= \int_{{\mathfrak g}^*/G} dP({\mathcal O})
\int_{\mathcal O} (\int_{\mathfrak g}
e^{i(\xi,X)}J^{-1/2}_{\mathfrak g}(X)\Phi(X)dX)
d\beta_{\mathcal O}(\xi),$$ where $P({\mathcal O})$ is the Plancher\`el
measure on ${\mathfrak g}^*/G$ and $\beta_{\mathcal O}(\xi)$ is the
Liouville measure on the orbit ${\mathcal O}$.
The proof of the main theorem is therefore achieved.

Let us now apply the construction of Chern-Weil homomorphism to our case
of $M^{\mathcal O}_0$. Let ${\mathcal A}(M^{\mathcal O}_0)$ be the algebra
of differential forms over $B = M^{\mathcal O}_0$. There is the natural
principal bundle $G \longrightarrow P^{\mathcal O} \longrightarrow
M^{\mathcal O}_0$, with connection 1-form $\omega = \sigma_{\mathfrak
g}^{\mathcal O}$ with curvature $\Omega = d\omega +
\frac{1}{2}[\omega,\omega] \in {\mathcal A}^2(P^{\mathcal O},{\mathfrak
g})^G$. The decomposition $TP^{\mathcal O} = T_hP^{\mathcal O} +
T_vP^{\mathcal O}$ of the tangent bundle $TP^{\mathcal O}$ into the sum of
horizontal and vertical components, $T_hP^{\mathcal O} = \ker(\omega)$.
The vertical sub-bundle $T_vP^{\mathcal O}$ determines a projection $h$
from ${\mathcal A}(P^{\mathcal O})$ onto the sub-algebra of horizontal
forms
$${\mathcal A}(P^{\mathcal O})_{hor} = \{ \alpha \in {\mathcal
A}(P^{\mathcal O}) \quad \vert \quad \imath(X)\alpha = 0, \forall X \in
{\mathfrak g} \}.$$ With respect to a basis $E_1, \dots, E_n$ of
${\mathfrak g}$, we can write 
$$\omega = \sum_{i=1}^n \omega^iE_i \quad \mbox{\rm and} \quad \Omega =
\sum_{i=1}^n\Omega^iE_i,$$ where $\omega^i \in {\mathcal A}^1(P^{\mathcal
O})_{hor}$ and $\Omega^i \in {\mathcal A}^2(P^{\mathcal O})_{hor}$ and the
projector $h$ can be written as
$$h = \prod_{i=1}^n (I-\omega^i \imath(E_i) = \sum_{1 \leq i_1 < \dots <
i_r} \omega^{i_1}\dots \omega^{i_r} \imath(E_{i_1})\dots
\imath(E_{i_r}),$$ see for example (\cite{berlinegetzlervergne},
Lemma 7.30). If ${\mathcal V} = {\mathcal E}_{\rho}(V) = P^{\mathcal O}
\times_G
V$ is an
vector bundle, associated with a representation $\rho$, with the
associated
affine connection, or equivalently,
a covariant derivation $\nabla^{\mathcal V}$, then $$\nabla^{\mathcal
V}\alpha = d\alpha + \sum_{i=1}^n \omega^i \rho(E_i)\alpha.$$

It is easy to see that the basic differential forms over $P^{\mathcal O}$
are just the ordinary differential forms on $B = M^{\mathcal O}_0 
= G\setminus P^{\mathcal O}$, ${\mathcal A}(P^{\mathcal O}) \cong
{\mathcal A}(M^{\mathcal O}_0)$,
$${\mathcal A}(P^{\mathcal O})_{bas} \cong {\mathcal A}(M^{\mathcal
O}_0);$$ the same for ${\mathcal V}$-valued
differential forms   
$${\mathcal A}(P^{\mathcal O},{\mathcal V}) = {\mathcal A}(P^{\mathcal
O},V)_{bas} \cong {\mathcal A}(M^{\mathcal
O}_0)\otimes_{\mathbb C}V.$$
If $\alpha \in {\mathcal A}(P^{\mathcal O},{\mathcal V})$ then $h\alpha =
\alpha$. This means that $h$ acts as the identity operator on 
${\mathcal A}(P^{\mathcal O},{\mathcal V}) = {\mathcal A}(P^{\mathcal O}, 
V)_{bas}$ and $d$ conserves this subspace of forms. We have
$$\imath(X)d\alpha = (\Lie(X) \otimes I)\alpha = -\rho(X)\alpha,$$ then
$$\imath(E_i)\imath(E_j)d\alpha \equiv 0, \forall E_i, E_j.$$ Thus,
$$h(dh\alpha) = h(d\alpha) = \nabla^{\mathcal V}\alpha.$$ It is easy to
see (\cite{berlinegetzlervergne}, Proposition 7.32) that the covariant
derivation $\nabla^{\mathcal V}$ on ${\mathcal V}$ coincides with the
restriction of $D \otimes I$ on ${\mathcal A}(P^{\mathcal O},V)$ to the
subspace ${\mathcal A}(P^{\mathcal O},{\mathcal V}) := {\mathcal
A}(P^{\mathcal O},V)_{bas}$, where $$D := hdh = h(d -\sum_{i=1}^n \Omega^i
\imath(E_i))$$ 
on ${\mathcal A}(P^{\mathcal O})$.

If $M$ is a $G$-manifold, denote ${\mathcal M} = P^{\mathcal O} \times_G
M$ the associated fibration. One has
$${\mathcal A}({\mathcal M}) \cong {\mathcal A}(P^{\mathcal O} \times_G
M)_{bas} =$$ $$= \{ \alpha \in {\mathcal A}(P^{\mathcal O} \times M)\quad \vert
\quad \alpha
\mbox{ is basic w.r.t. the action of G} \},$$
$$\imath(X) := \imath(X_{P^{\mathcal O} \times M}) = \imath(X_{P^{\mathcal
O}}) + \imath(X_M),$$
$h : {\mathcal A}(P^{\mathcal O} \times M) \to {\mathcal A}(P^{\mathcal
O} \times M)_{bas}$ is given by
$$h = \prod_{i=1}^m (I - \omega^i \imath(X_i)),$$
$$D = h.d.h|_{{\mathcal A}({\mathcal M})} = h.d.h|_{{\mathcal
A}(P^{\mathcal O} \times M)_{bas}} = d_{\mathcal M}.$$
Let us consider the complex $({\mathbb C}[{\mathfrak g}] \otimes
{\mathcal A}(M), d_{\mathfrak g})$. We correspond to each element $\alpha
= f \otimes \beta$ the element $\alpha(\Omega) := f(\Omega) \otimes \beta
\in {\mathcal A}(P^{\mathcal O}) \otimes {\mathcal A}(M)$. One defines
thus the so called {\it Chern-Weil homomorphism} 
$$W=\Phi_\omega : {\mathbb C}[{\mathfrak g}] \otimes {\mathcal A}(M) \to
{\mathcal A}(P^{\mathcal O} \times M)_{hor},$$
$$W= \Phi_\omega(\alpha) := h(\alpha(\Omega)).$$ It was shown
(\cite{berlinegetzlervergne},Theorem 7.34) that if $G$ is a Lie group, $G
\longrightarrow P^{\mathcal O} \longrightarrow M^{\mathcal O}_0$ is the
principal bundle with connection $\omega = \sigma_{red}^{\mathcal O}$ and
the curvature $\Omega = d\omega + \frac{1}{2}[\omega,\omega]$ and $M$ is a
smooth $G$-manifold, then the Chern-Weil homomorphism induces a
homomorphism  of the complexes of differential graded algebras
$$({\mathcal A}_G(M), d_{\mathfrak g}) = (({\mathbb C}[{\mathfrak g}]
\otimes {\mathcal A}(M))^G,d_{\mathfrak g}) \to ({\mathcal
A}({\mathcal M}),d).$$ In particular, if $M = \{pt\}$ is a point, we obtain
the classical Chern-Weil homomorphism
$$W=\Phi_\omega : {\mathbb C}[{\mathfrak g}]^G \to {\mathcal
A}(M^{\mathcal
O}_0).$$ If $M=V$ is a vector $G$-space, we have Chern-Weil homomorphism
for induced vector bundle
$$W=\Phi_\omega : f \otimes v \mapsto f(\Omega) \otimes v,$$ which is a
homomorphism from $G$-equivariant K-groups to the even de Rham cohomologie
$\oplus_* H^{2*}_{DR}(M^{\mathcal O}_0; {\mathbb R})$. 

We can now reformulate the main result as what follows.
\begin{thm} If $\Phi$ is an $G$-invariant, then
$$\int_{\mathfrak g}\Theta_0^{\mathcal O}(X)\Phi(X)dX = (2\pi
i)^{\dim G}\vol(G)\int_{M^{\mathcal O}_{red}}\alpha_{red}W(\Phi),$$ where
$W : C^{\infty}({\mathfrak g})^G \to H^*(M^{\mathcal O}_{red})$ is the
Chern-Weil homomorphism, associated to the principal fibration
$\mu^{-1}({\mathcal O}) \to M^{\mathcal O}_{red}$.  
\end{thm}

\section*{Acknowledgments} The author would like to thank M. Vergne for
the related reprints of her works she sent. 

This work is supported in part by the International Centre for Theoretical
Physics (Italy), and the Vietnam National Research Programme in
Fundamental Sciences. The final version of this work was completed
during a stay of the author at the UP Diliman, the Philippines. The author
would like to thank Department of Mathematics and UP Diliman for
hospitality and especially professor Milagros P. Navarro for the
invitation, help and friendship.

\bibliographystyle{amsplain}
\bibliography{References}

\begin{thebibliography}{lbl}

\bibitem{berlinegetzlervergne} {\sc N. Berline, E. getzler and Mich\`ele
Vergne} \textit{Heat Kernels and Dirac Operators}, Grundlehren der mathematischen
Wischenschaften, No 298, Springer-Verlag, 1994; 2nd correcting printing 1996.

\bibitem {guilleminsternberg} {\sc V. Guillemin} and {\sc S. Sternberg},
\textit{Symplectic Techniques in Physics}, Cambridge University Press,
Cambridge, 1984. 

\bibitem {jeffreykirwan} {\sc L. C. Jeffrey} and {\sc F. C. Kirwan},
\textit{ Localization for nonabelian group action}, Topology
\textbf{34}(1995) 291-268. 

\bibitem {kalkman} {\sc J. Kalkman}, \textit{Cohomology rings of
symplectic quotients}, preprint, University Utrecht 1993. 

\bibitem{kirillov}
{\sc A. A. Kirillov}, \textit{Elements of the Theory of Representations},
Springer-Verlag, Berlin-New York,-Heidelberg, 1975.

\bibitem {vergne} {\sc M. Vergne}, \textit{A note on the
Jeffrey-Kirwan-Witten localisation formula}, Topology, \textbf{35}(1996),
No 1, 243-266. 

\bibitem {wu} {\sc S. Wu}, \textit{An integration formula for the square
of moment maps of circle actions}, preprint Hep-th/921207. 

\end{thebibliography}

\end{document}